\documentclass[11pt]{amsart}
\usepackage[margin=1.1in]{geometry}
\usepackage{amsmath,amssymb,amsthm,mathtools}
\usepackage{booktabs}
\usepackage{enumitem}
\usepackage{xcolor}
\usepackage{breqn}
\usepackage[colorlinks=true,linkcolor=blue!50!black,citecolor=blue!50!black,urlcolor=blue!50!black]{hyperref}

\newtheorem{theorem}{Theorem}[section]
\newtheorem{proposition}[theorem]{Proposition}
\newtheorem{lemma}[theorem]{Lemma}
\newtheorem{corollary}[theorem]{Corollary}
\newtheorem*{mainthm}{Theorem}
\newtheorem*{maincor}{Corollary}
\theoremstyle{definition}
\newtheorem{definition}[theorem]{Definition}
\theoremstyle{remark}
\newtheorem{remark}[theorem]{Remark}
\newtheorem{example}[theorem]{Example}

\newcommand{\F}{\mathbb{F}}
\newcommand{\Q}{\mathbb{Q}}
\newcommand{\Z}{\mathbb{Z}}
\newcommand{\PP}{\mathbb{P}}
\newcommand{\LL}{\mathcal{L}}
\newcommand{\gon}{\operatorname{gon}}
\newcommand{\divv}{\operatorname{div}}
\newcommand{\supp}{\operatorname{supp}}
\newcommand{\Aut}{\operatorname{Aut}}

\title{The $\Q$-gonality of $X_1(42)$ is $12$}
\author{Sourajyoti Maiti}
\address{Department of Mathematics, University of Connecticut, Storrs, CT 06269, USA}
\email{sourajyoti.maiti@uconn.edu}
\thanks{This work was partially supported by US National Science Foundation,  grant DMS-2617278.}
\date{\today}
\begin{document}
\begin{abstract}
The $\Q$-gonality of the modular curve $X_1(N)$ is known for all $N\le 40$
by work of Derickx and van Hoeij. We prove that the $\F_5$-gonality and the
$\Q$-gonality of $X_1(42)$ are both $12$. The upper bound $12$ was already known, and we raise the previous
lower bound from $10$ to $12$. Since the set of rational points of the Jacobian, $J_1(42)(\Q)$, is finite, it follows that
$12$ is the least degree in which $X_1(42)$ has infinitely many points. The
lower bound comes from an exhaustive search over Riemann--Roch spaces on
$X_1(42)$ over $\F_5$, following Derickx and van Hoeij. We reduce this search
by the orbits of a cyclic group of automorphisms of order $12$ on places of
degree at most $5$.  This leaves $743\,231$ Riemann--Roch spaces, about a quarter of the
number needed by the earlier method. We prove the reductions and provide the
code and data needed to reproduce the computations.
\end{abstract}

\maketitle

\section{Introduction}

Let $k$ be a field and $X/k$ a smooth, projective, geometrically irreducible
curve. The \emph{$k$-gonality} $\gon_k(X)$ is the least degree of a nonconstant
morphism $X\to\PP^1_k$, i.e.\ the least degree of a nonconstant $f\in k(X)$.
For the modular curves $X_1(N)$ the $\Q$-gonality controls the arithmetic of
points of small degree. By a theorem of Frey \cite[Proposition~2]{Frey}, if $K$ is a number field, $X/K$ has a $K$-rational point and $X$ has
infinitely many points of degree $\le d$ over $K$, then $\gon_K(X)\le 2d$.
Lower bounds for $\gon_\Q X_1(N)$ therefore give finiteness results: if
$\gon_\Q X_1(N)>2d$, then there are only finitely many elliptic curves (up to
$\overline{\Q}$-isomorphism) with a point of order $N$ over number fields of
degree $\le d$. Such bounds are a standard
ingredient in classifications of torsion groups over number fields of small
degree; see e.g.~\cite{DvH,Najman}.

Derickx and van Hoeij~\cite{DvH} developed a method for computing
$\gon_{\F_p}X_1(N)$ by an exhaustive search over Riemann--Roch spaces, and
determined $\gon_\Q X_1(N)$ exactly for all $N\le 40$; the cases $N\le 22$
were known before~\cite[p.~2]{SuthMC}. For $N\le 250$ they gave upper
bounds by explicit functions over $\Q$, found by an LLL search among the
modular units \cite[Theorem~4]{DvH}; for $N=42$ the bound is
\begin{equation}\label{eq:upper}
  \gon_\Q X_1(42)\le 12
\end{equation}
\cite[Table~1]{DvH}. On the side of lower bounds, Derickx and van Hoeij's
computation~\cite{vHX142}, which accompanies~\cite{DvH}, shows that $X_1(42)_{\F_5}$
has no nonconstant function of degree $\le 8$, so $\gon_\Q X_1(42)\ge 9$. This
is the bound $\gon_\Q X_1(42)>8$ used in \cite[Proof of Proposition~6.3]{Najman}.
Using the new algorithm of Derickx and Terao~\cite{DT}, Najman and Varivoda
improved it to $\gon_{\F_5}X_1(42)\ge 10$, and hence $\gon_\Q X_1(42)\ge 10$
\cite[Table~1]{NV}; this bound is also quoted in~\cite{Isolated}. To our
knowledge this was the best known lower bound.

\par Lower bounds for the $\Q$-gonality are obtained from the following
inequality. If $C/\Q$ is a curve and $p$ is a prime of good reduction for $C$,
then
\begin{equation}\label{eq:red}
  \gon_{\F_p}(C_{\F_p})\ \le\ \gon_\Q(C).
\end{equation}
In this form, \eqref{eq:red} was proved by Derickx in his master's
thesis \cite[Theorem~2.5]{DerickxMSc}, where it is stated for a projective
scheme over a discrete valuation ring that is smooth of relative dimension one
with geometrically irreducible fibres; see also \cite[Section~1, (1)]{DvH}. A similar statement
was given earlier by Frey~\cite{Frey}, who attributes it to Deuring's work on
the reduction of algebraic function fields~\cite{Deuring}. The curve $X_1(N)$ has good
reduction at every prime $p\nmid N$ \cite[Section~1]{DvH}.

\par The difficulty in improving the gonality lower bounds was the feasibility of computation. For example, the earlier method of Derickx and van Hoeij \cite{DvH,vHX142}, run on a single CPU, would take more than four months (in wall-clock time) to decide whether $X_1(42)_{\F_5}$ has a nonconstant function of degree $\le 11$, and it also runs into memory limits. Our computation is built on their method, but we use different techniques to reduce the search and make the computation feasible (our best search technique took about 1.9 days in wall-clock time to decide the same). The main contribution of this paper is the exact determination of the gonality; the orbit-reduced search of Section~\ref{sec:orbitsearch} is a computational refinement of the method of~\cite{DvH}. In particular, $N=42$ is
the first $N>40$ for which $\gon_\Q X_1(N)$ is now known. 
 Our main result is the following.

\begin{mainthm}
Every nonconstant function in $\F_5(X_1(42))$ has degree at least $12$. Consequently
\[
  \gon_{\F_5}X_1(42)\;=\;\gon_{\Q}X_1(42)\;=\;12 .
\]
\end{mainthm}

\begin{maincor}
$X_1(42)$ has only finitely many points of degree $\le 11$, and
infinitely many points of degree $12$. Thus $12$ is the least $d$ such
that $X_1(42)$ has infinitely many points of degree $d$. In particular,
only finitely many $\overline{\Q}$-isomorphism classes of elliptic curves over
number fields of degree $\le 11$ have a point of order $42$.
\end{maincor}

The finiteness of the points of degree $\le 5$ follows from Frey's criterion
and was already known \cite[Section~5, Lemma~1 and Proposition~6]{DvH}. The stronger statement uses, in addition, that the Mordell--Weil group
$J_1(42)(\Q)$ has rank $0$ \cite[Theorem~3.1(2)]{DEvMZB}; see also
\cite[Theorem~2.3(4)]{NV}.

\begin{proof}[Proof of the corollary]
Let $C$ be a nice curve over a number field $K$, that is, a smooth,
projective, geometrically integral curve over $K$. If the Jacobian of $C$
has rank $0$ over $K$, then the least $d$ such that $C$ has infinitely many
points of degree $d$ is $\gon_K C$ \cite[Proposition~2.3]{DS}. Here points of
degree $d$ are points $a\in C(\overline{K})$ with $[K(a):K]=d$. Since such
points come in Galois orbits of size $d$, the statement is the same whether
they are counted individually or up to Galois conjugacy.

The curve $X_1(42)$ is nice over $\Q$, its Jacobian has rank $0$ over $\Q$
\cite[Theorem~3.1(2)]{DEvMZB}, and $\gon_\Q X_1(42)=12$ by the Theorem.
Hence $12$ is the least degree in which $X_1(42)$ has infinitely many
points: for every $d\le 11$ there are only finitely many points of
degree $d$, and there are infinitely many of degree $12$.

For the last statement, a pair $(E,P)$ with $E$ defined over a number field
$L$ of degree $\le 11$ and $P\in E(L)$ of order $42$ gives a non-cuspidal
point of $X_1(42)$ of degree $\le 11$, so only finitely many such points
occur. Since the forgetful map $X_1(42)\to X(1)\cong\PP^1_j$,
$(E,P)\mapsto j(E)$, assigns one $j$-invariant to each point, and elliptic
curves over $\overline{\Q}$ are classified up to isomorphism by their
$j$-invariant, only finitely many $\overline{\Q}$-isomorphism classes of $E$
occur. 
\end{proof}

\subsection*{Method}
We follow \cite{DvH}. Let $S$ be the sum of the $n=12$ rational points of
$X_{\F_5}=X_1(42)_{\F_5}$. A M\"obius transformation over $\F_5$ moves any
function of degree $\le d$ to one whose pole divisor contains at least
$\lceil n/6\rceil=2$ rational points (Lemma~\ref{lem:mobius}). Every such
function lies in $\LL(S+E)$ for an effective divisor $E$ of degree exactly
$d-2$ (Lemma~\ref{lem:cover}). We use the automorphisms of $X_{\F_5}$ to reduce the number of divisors to check (Lemma~\ref{lem:aut}). We enumerate these divisors $E$ by
the partition of $d-2$ formed by the degrees of their places, and compute
$\LL(S+E)$ for each with Hess's algorithm~\cite{Hess} in Magma~\cite{Magma}.

\par We do this in stages. As we only knew the lower bound from \cite{vHX142} at the time we started this project, we started with $d=10$ and tested $600\,465$ Riemann--Roch spaces
and found no nonconstant function of degree $\le 10$. We then ran the search for
$d=11$, testing $2\,994\,534$ Riemann--Roch spaces, and found no nonconstant
function of degree $11$ (nor of any smaller degree). These searches use the automorphisms only through their action on the rational points, a reduction taken directly from Derickx and van Hoeij's code~\cite{vHX142}. We then reduced the search for $d=11$ further, to $743\,231$ divisors, by using a cyclic group of automorphisms of order $12$, found with Magma, and its action on places of higher degree. This orbit-reduced search is the one on which our proof rests: the earlier searches split the work among jobs by the positions of places in lists computed by Magma, whose order is not reproducible between sessions, so for them we cannot certify that every divisor was tested, and we regard them as supporting evidence. All these searches were split into $20$ independent jobs
run in parallel. The code was written with the assistance of AI and reviewed
in full by the author (see the \hyperref[sec:ai]{Declaration of generative AI use}).

The upper bound requires no new computation. It is the bound
\eqref{eq:upper} of Derickx and van Hoeij \cite[Table~1]{DvH}, given by an
explicit modular unit $g_{42}$ of degree $12$ over $\Q$, which is listed on
the authors' webpage. Derickx and van Hoeij note that their probabilistic LLL
search cannot prove the optimality of the functions $g_N$, but that the bounds
in their Table~1 are likely sharp whenever the $\Q$-gonality is attained by a
modular unit \cite[Section~1]{DvH}; for $N=42$ our theorem confirms that the
bound is sharp. The same bound can also be read off the optimized plane model of $X_1(42)$ in
Sutherland's tables~\cite{SuthOpt}.
Together with \eqref{eq:red} this gives
$\gon_{\F_5}X_1(42)\le\gon_\Q X_1(42)\le 12$, so our lower bound is sharp over
both fields.

In a separate project, we tried to improve Sutherland's optimized model for $X_1(42)$. We used several search techniques, including reinforcement-learning-guided search, which did not succeed. These heuristics suggested that the gonality should be exactly $12$ and motivated the exhaustive search carried out here.

\subsection*{Outline}
Section~\ref{sec:prelim} collects facts about $X_1(42)$ over $\F_5$.
Section~\ref{sec:criterion} proves the search criterion and
Section~\ref{sec:sym} the basic lemma on automorphisms. Section~\ref{sec:tech}
describes the first search, which serves as supporting evidence. Section~\ref{sec:orbitsearch} is the orbit-reduced search on which the proof rests, and
Section~\ref{sec:proof} completes the proof.
\section{Preliminaries}\label{sec:prelim}

\subsection{Gonality and reduction}
We use standard notation for a curve $X$ over a field $k$: divisors are
finite sums $D=\sum_Q n_QQ$ of places $Q$ of $k(X)$, with
$\deg D=\sum_Q n_Q[\kappa(Q):k]$, where $\kappa(Q)$ is the residue field of
$Q$; $(f)_0$ and $(f)_\infty$ are the zero and pole divisors of
$f\in k(X)^\times$, so that $\divv f=(f)_0-(f)_\infty$; for nonconstant $f$,
$\deg f=\deg(f)_\infty=[k(X):k(f)]$ is the degree of $f\colon X\to\PP^1_k$;
and for $D\ge0$, $\LL(D)=\{f\in k(X)^\times:\divv f+D\ge0\}\cup\{0\}$ and
$\ell(D)=\dim_k\LL(D)$. A nonconstant $f$ lies in $\LL(D)$ if and only if
$(f)_\infty\le D$, and then $\deg f\le\deg D$.

Inequality~\eqref{eq:red} is the basis of the method of~\cite{DvH}. To
determine $\gon_\Q X_1(N)$ one needs a function over $\Q$ whose degree is an
upper bound, together with a prime $p\nmid N$ for which an exhaustive search
shows that $\gon_{\F_p}X_1(N)$ equals that degree; see the discussion of
Tasks~1--3 in \cite[Section~1]{DvH}. Derickx and van Hoeij observed that for
every $N\le 40$ the $\Q$-gonality equals the $\F_p$-gonality for the smallest
prime $p\nmid N$ \cite[Remark~1]{DvH}, although this cannot be expected for
all $p$. Our result shows that the same holds for $N=42$, where the smallest
such prime is $p=5$.

We apply \eqref{eq:red} with $C=X_1(42)$ and $p=5$, and write
$\boxed{X=X_1(42)_{\F_5}}$ from now on.

\subsection{Automorphisms}
Let $C$ be a curve over a field $k$.

\begin{definition}\label{def:aut}
The \emph{automorphism group of $C$ over $k$} is
\[
  \Aut_k(C)=\bigl\{\sigma\colon C\to C \;:\; \sigma \text{ is an isomorphism of } k\text{-schemes}\bigr\},
\]
with composition as group law. The group of $k$-automorphisms of the function
field is
\[
  \Aut\bigl(k(C)/k\bigr)=\bigl\{\tau\colon k(C)\to k(C) \;:\; \tau \text{ is a field automorphism with } \tau|_k=\mathrm{id}\bigr\}.
\]
\end{definition}

Every $\sigma\in\Aut_k(C)$ induces $\sigma^*\in\Aut(k(C)/k)$, $\sigma^*(f)=f\circ\sigma$,
and since $C$ is smooth and projective this gives an anti-isomorphism
\[
  \Aut_k(C)\xrightarrow{\ \sim\ }\Aut\bigl(k(C)/k\bigr),\qquad \sigma\longmapsto\sigma^*,
  \qquad (\sigma\circ\rho)^*=\rho^*\circ\sigma^* .
\]
For $P\in C$,\, $\sigma(P)$ is the
place with valuation
\[
  v_{\sigma(P)}(f)=v_P(\sigma^*f)\qquad\text{for all } f\in k(C)^\times .
\]
This action of $\sigma$ preserves the degree of places, in particular it permutes $C(k)$. It
extends linearly to divisors,
\[
  \sigma\Bigl(\sum_Q n_QQ\Bigr)=\sum_Q n_Q\,\sigma(Q),
\]
preserving degree and effectivity, and $\divv(f\circ\sigma^{-1})=\sigma(\divv f)$
for $f\in k(C)^\times$. If $g(C)\ge2$, then $\Aut_k(C)$ is finite
\cite[Theorem~1.11]{DM}.
\subsection{Modular curves and diamond operators}
Let $N\ge 4$. The affine modular curve $Y_1(N)$ over $\Q$ is the fine moduli
space of pairs $(E,P)$, where $E$ is an elliptic curve and $P\in E$ is a point
of exact order $N$, two pairs $(E,P)$ and $(E',P')$ being identified if there
is an isomorphism $E\to E'$ taking $P$ to $P'$. Its smooth projective
completion is $X_1(N)$, and the points of $X_1(N)\setminus Y_1(N)$ are the
\emph{cusps}. The curve $X_1(N)$ has good reduction at every prime $p\nmid N$.
For such $p$, the curve $Y_1(N)_{\F_p}$ has the same moduli interpretation
over fields of characteristic $p$, and $X_1(N)_{\F_p}$ is its smooth
projective completion; its cusps do not correspond to pairs $(E,P)$ with $E$
an elliptic curve.

\begin{definition}\label{def:diamond}
For $a\in(\Z/N\Z)^\times$, the \emph{diamond operator} $\langle a\rangle$ is
the automorphism of $X_1(N)$ given on $Y_1(N)$ by
\[
  \langle a\rangle\colon (E,P)\longmapsto (E,aP)
\]
and extended to the cusps. The \emph{diamond group} is
\[
  G_N=\bigl\{\langle a\rangle \;:\; a\in(\Z/N\Z)^\times\bigr\}.
\]
\end{definition}

The diamond operators are defined over $\Q$, and $\langle ab\rangle=\langle
a\rangle\circ\langle b\rangle$. Since $[-1]\colon E\to E$ is an isomorphism
taking $P$ to $-P$, we have $\langle -1\rangle=\mathrm{id}$, so
\[
  G_N\cong(\Z/N\Z)^\times/\{\pm1\}\ \subseteq\ \Aut_\Q(X_1(N)),
\]
and the quotient $X_1(N)/G_N$ is $X_0(N)$. For $p\nmid N$ the diamond
operators reduce to automorphisms of $X_1(N)_{\F_p}$ defined over $\F_p$, with
the same moduli description.

\subsection{Basic invariants}
The genus of $X_1(42)$ is
\[
  g=1+\frac{42^2}{24}\prod_{p\mid42}\Bigl(1-\frac1{p^2}\Bigr)-\frac14\sum_{e\mid42}\varphi(e)\varphi(42/e)
   =1+48-24=25,
\]
and $X_1(42)$ has $48$ cusps. Its diamond group is
$G_{42}\cong(\Z/42\Z)^\times/\{\pm1\}\cong\Z/6\Z$
(Definition~\ref{def:diamond}), and since $5\nmid 42$ it acts on $X$ by
automorphisms defined over $\F_5$.

\begin{lemma}\label{lem:cusps}
Every point of $X(\F_5)$, and more generally of $X(\F_{25})$, is a cusp.
\end{lemma}

\begin{proof}
A non-cuspidal point of $X$ over $\F_q$ gives an elliptic curve $E/\F_q$ with
a point of order $42$ in $E(\F_q)$, since $Y_1(42)$ is a fine moduli space.
The Hasse bound gives $\#E(\F_5)\le 5+1+2\sqrt5<11$ and
$\#E(\F_{25})\le 25+1+10=36<42$.
\end{proof}

\subsection{The model and the places of small degree}
We use the plane model $f(x,y)=0$ of $X_1(42)$ in Sutherland's
tables~\cite{Sutherland} (file \texttt{FFFc42}), where $f\in\Z[x,y]$. Let
$\bar f\in\F_5[x,y]$ be its reduction modulo $5$; this is the polynomial used
in all our computations. Let $F$ be the function field of the plane curve
$\bar f=0$, formed in Magma as an algorithmic function field. All
computations take place in $F$, so the singularities of the plane model play
no role. That $F$ is the function field of $X$ is not automatic, since the
reduction of a plane model need not be the reduction of the curve. For example, the elliptic curve $y^2=x^3+1$ over $\Q$ has good
reduction at $5$, and the substitution $x=U/25$, $y=V/125$ gives the
birationally equivalent model $V^2=U^3+5^6$ with integer coefficients.
Modulo $5$ this model becomes $V^2=U^3$, which is absolutely irreducible but
has genus $0$: its function field is $\F_5(V/U)$, not the function field of
the reduction of the elliptic curve. We prove that $F$ is the function field of $X$ by
using the moduli interpretation of $X_1(42)$ over $\F_5$: Lemma~\ref{lem:criterion}
gives a criterion for a function field over $\F_5$ to be that of $X$, and
Lemma~\ref{lem:model} verifies it for $F$. Lemma~\ref{lem:criterion} is standard in substance. We include a proof of the form we need.

For $b,c$ in a field, let $E(b,c)$ denote the curve in Tate normal form
\[
  E(b,c)\colon\ V^2+(1-c)UV-bV=U^3-bU^2,
\]
with the point $(U,V)=(0,0)$.

\begin{lemma}\label{lem:criterion}
Let $L$ be a function field of one variable over $\F_5$. Suppose that there
are $b,c\in L$ such that
\begin{enumerate}[label=(\alph*),leftmargin=*]
  \item $E(b,c)$ is an elliptic curve over $L$ on which $(0,0)$ has exact
    order $42$, and
  \item $L=\F_5(b,c)$.
\end{enumerate}
Then $L\cong\F_5(X)$.
\end{lemma}

\begin{proof}
Let $C$ be the smooth projective curve over $\F_5$ with function field $L$.

Recall the Tate normal form: if $L'$ is a field and $P'\in E'(L')$ is
a point of order at least $4$ on an elliptic curve $E'/L'$, then there are
unique $b',c'\in L'$ such that $(E',P')$ is isomorphic over $L'$ to
$\bigl(E(b',c'),(0,0)\bigr)$ \cite{Kubert,Husemoller}. By uniqueness,
$(b',c')$ does not change when $(E',P')$ is base changed to a larger field.

Since $5\nmid42$, the curve
$Y_1(42)_{\F_5}$ is a fine moduli scheme for pairs $(E',P')$ consisting of an
elliptic curve $E'$ and a point $P'$ of exact order $42$~\cite{KM}, and $X$ is
the smooth, projective, geometrically irreducible curve containing
$Y_1(42)_{\F_5}$ as a dense open subset (\cite{DR,KM}; see
\cite[Theorems~3.2.7 and 4.2.1]{Conrad}). Put $E=E(b,c)$ and $P=(0,0)$. By
(a), the pair $(E,P)$ defines a morphism $\operatorname{Spec}L\to
Y_1(42)_{\F_5}$, that is, a rational map $C\dashrightarrow X$, which extends
to a morphism $\varphi\colon C\to X$ because $C$ is a smooth curve and $X$ is
proper.

 The morphism $\varphi$ is nonconstant. Indeed, by (b), $b$ and $c$ are not both
algebraic over $\F_5$, since $L$ has transcendence degree $1$ over $\F_5$. If $\varphi$ were constant, with image a closed point $x_0\in X$
with (finite) residue field $\kappa(x_0)$, then $(E,P)$ would be the base
change to $L$ of a pair defined over $\kappa(x_0)$, and by the uniqueness of the Tate normal form we would
have $b,c\in\kappa(x_0)$, a contradiction. Hence $\varphi$ maps the generic
point of $C$ to the generic point $\eta$ of $X$ and induces an embedding of
fields $\varphi^*\colon\F_5(X)\hookrightarrow L$.

We now show that $b$ and $c$ lie in the image of $\varphi^*$. The generic point
$\eta$ corresponds to a pair $(E_\eta,P_\eta)$ over $\F_5(X)$ with $P_\eta$
of exact order $42$, so it has a Tate normal form with parameters
$b_\eta,c_\eta\in\F_5(X)$. Since $Y_1(42)_{\F_5}$ is a fine moduli scheme, the
base change of $(E_\eta,P_\eta)$ along $\varphi^*$ is isomorphic to $(E,P)$,
and by the uniqueness of the Tate normal form, $\varphi^*(b_\eta)=b$ and
$\varphi^*(c_\eta)=c$.

Consequently, the image of $\varphi^*$ is a subfield of $L$
containing $b$ and $c$, which generate $L$ by (b). Hence $\varphi^*$ is
surjective, and therefore an isomorphism: $L\cong\F_5(X)$.
\end{proof}

\begin{lemma}\label{lem:model}
$F\cong\F_5(X)$, that is, the normalization of the curve $\bar f=0$ is
$X=X_1(42)_{\F_5}$.
\end{lemma}

\begin{proof}
Let $x,y\in F$ be the coordinate functions, so that $\bar f(x,y)=0$.
Sutherland~\cite{Sutherland} gives a birational map from $f(x,y)=0$ to the
curve defined by the Tate normal form, together with its inverse. Following
him, put
\[
  r=\frac{x^2y-xy+y-1}{x^2y-x},\qquad s=\frac{xy-y+1}{xy},\qquad
  b=rs(r-1),\qquad c=s(r-1).
\]
With Magma we verified the following (script \texttt{X1\_42\_moduli\_check.m}
in the data repository):
\begin{enumerate}[label=(\roman*),leftmargin=*]
  \item $\bar f$ is irreducible over $\F_5$, so $F$ is a function field of
    one variable over $\F_5$;
  \item the denominators of $r$ and $s$ are nonzero in $F$, the discriminant
    of $E(b,c)$ is nonzero, so $E(b,c)$ is an elliptic curve over $F$, and for
    $P=(0,0)$ we have $[42]P=O$ while $[21]P$, $[14]P$ and $[6]P$ are
    nonzero, so $P$ has exact order $42$;
  \item Sutherland's inverse map recovers the coordinates:
    \[
      x=\frac{s-r}{rs-2r+1},\qquad y=\frac{rs-2r+1}{s^2-s-r+1}\qquad
      \text{in } F.
    \]
\end{enumerate}
By (ii), hypothesis (a) of Lemma~\ref{lem:criterion} holds for $L=F$. By
(iii), $x,y\in\F_5(r,s)$, so $F=\F_5(x,y)=\F_5(r,s)$. Since $P=(0,0)$ has exact order $42$, we have $c\ne0$ and $b\ne c$: on $E(b,c)$
one has $-2P=(b,0)$, $2P=(b,bc)$ and $3P=(c,b-c)$, so $c=0$ would give $4P=O$
and $b=c$ would give $5P=O$. Hence $r=b/c$ and $s=c^2/(b-c)$. Hence $F=\F_5(b,c)$, which is hypothesis (b). By
Lemma~\ref{lem:criterion}, $F\cong\F_5(X)$.
\end{proof}

From now on we identify $F$ with $\F_5(X)$.

Let $P_v$ be the number of places of $F$ of degree $v$. Magma gives
\begin{center}
\begin{tabular}{@{}c|ccccccccc@{}}
\toprule
$v$ & 1 & 2 & 3 & 4 & 5 & 6 & 7 & 8 & 9\\
\midrule
$P_v$ & 12 & 6 & 8 & 183 & 648 & 2288 & 11160 & 48708 & 214384\\
\bottomrule
\end{tabular}
\end{center}
In particular $n=\#X(\F_5)=12$. The values for $v\le6$ agree with those
recorded in~\cite{vHX142}, and those for $v\le 5$ with Najman's description
of $X_1(42)$ modulo $5$ \cite[Remark~6.4 and Table~3]{Najman}. There, the Galois
orbits of cusps over $\Q$ have degrees $1^{12},2^6,6^4$; modulo $5$ two of the
four orbits of degree $6$ split into two places of degree $3$ each. Table~3
of \cite{Najman} lists $0,0,4,183,648$ non-cuspidal places of degrees $1,\dots,5$.
Adding the cuspidal places $12,6,4,0,0$ gives $P_1,\dots,P_5$.

\section{The search criterion}\label{sec:criterion}
The results of this section systematically organize ideas from the work of
Derickx and van Hoeij, in particular \cite[Proposition~5]{DvH} and van
Hoeij's code~\cite{vHX142}.
Throughout, let
\[
  S=\sum_{Q\in X(\F_5)}Q
\]
be the sum of the $n=12$ rational places of $X$, each with multiplicity one.
By Lemma~\ref{lem:cusps} every rational place is a cusp, so $S$ is the sum of
the $\F_5$-rational cusps; following~\cite{DvH} we call it the
\emph{cusp sum} (\texttt{cuspsum} in the Magma scripts). It is an effective
divisor of degree $\deg S=12$. As before, $(h)_\infty$ denotes the pole
divisor of a nonconstant function $h$, so that $\deg h=\deg\,(h)_\infty$ and
$h\in\LL(D)\iff(h)_\infty\le D$ (Section~\ref{sec:prelim}).

For an effective divisor $D=\sum n_QQ$ put
$D\wedge S=\sum_{Q\in X(\F_5)}\min(n_Q,1)\,Q$, the sum of the distinct
rational points in $\supp D$. Then $D\wedge S\le D$ and $D\wedge S\le S$.

\begin{lemma}[M\"obius normalization]\label{lem:mobius}
Let $f\in\F_5(X)$ be nonconstant of degree $d'$. There is $h\in\F_5(X)$ with
$\deg h=d'$ and $\deg\bigl((h)_\infty\wedge S\bigr)\ge\lceil n/6\rceil=2$.
\end{lemma}

\begin{proof}
The morphism $f\colon X\to\PP^1$ maps the $12$ points of $X(\F_5)$ to the $6$
points of $\PP^1(\F_5)$, so by the pigeonhole principle some fibre
$f^{-1}(a)$, $a\in\PP^1(\F_5)$, contains at least $\lceil 12/6\rceil=2$
rational points. Put $h=1/(f-a)$ if $a\in\F_5$ and
$h=f$ if $a=\infty$. Then $\deg h=\deg f$, and every point of $f^{-1}(a)$ is a
pole of $h$.
\end{proof}

\begin{lemma}[Covering]\label{lem:cover}
Let $d\ge 2$, $B=d-2$, and let $h$ be nonconstant with $\deg h\le d$ and
$\deg((h)_\infty\wedge S)\ge 2$. There is an effective divisor $E$ with
$\deg E=B$ such that $h\in\LL(S+E)$.
\end{lemma}

\begin{proof}
Let $D=(h)_\infty$; this is an effective divisor of degree $\deg h\le d$.
Since $D\wedge S\le D$, the divisor $E'=D-(D\wedge S)$ is effective, and
$\deg E'=\deg D-\deg(D\wedge S)\le d-2=B$. Choose a rational place $Q_0$ and
put $E=E'+(B-\deg E')\,Q_0$. Then $E\ge E'\ge0$ and $\deg E=B$. Since
$D\wedge S\le S$ and $E'\le E$, we have $(h)_\infty=(D\wedge S)+E'\le S+E$.
As $(h)_0\ge0$, it follows that
$\divv h+S+E=(h)_0+\bigl(S+E-(h)_\infty\bigr)\ge0$, that is,
$h\in\LL(S+E)$.
\end{proof}

\begin{corollary}\label{cor:criterion}
Let $d\ge2$. Then $\gon_{\F_5}(X)\ge d+1$ if and only if, for every effective
divisor $E$ of degree $d-2$, every nonconstant function in $\LL(S+E)$ has
degree $>d$.
\end{corollary}

\begin{proof}
Suppose $\gon_{\F_5}(X)\le d$, and let $f$ be nonconstant of degree
$d'\le d$. By Lemma~\ref{lem:mobius} there is $h$ with $\deg h=d'$ and
$\deg((h)_\infty\wedge S)\ge2$, and by Lemma~\ref{lem:cover} there is an
effective $E$ of degree $d-2$ with $h\in\LL(S+E)$; so the condition fails.
Conversely, if the condition fails, some $\LL(S+E)$ contains a nonconstant
function of degree $\le d$, and $\gon_{\F_5}(X)\le d$.
\end{proof}

Note that $\deg(S+E)=n+d-2\le 21$ for $d\le 11$, which is less than $g$, so
Riemann--Roch does not force $\LL(S+E)$ to contain nonconstant functions; in
most cases $\ell(S+E)=1$.

\section{Reduction by automorphisms}\label{sec:sym}

We use Definition~\ref{def:aut} with $k=\F_5$. Since $g(X)=25$, the group
$\Aut_{\F_5}(X)$ is finite. Magma computes the elements of
$\Aut(\F_5(X)/\F_5)$; via $\sigma\mapsto\sigma^*$ we write $\sigma$ for both
when no confusion can arise. Concretely, Magma returns field automorphisms
$\sigma_1,\dots,\sigma_{12}$ of $\F_5(X)$ (Appendix~\ref{app:aut}); we write
$\varphi_i$ for the automorphism of $X$ with $\varphi_i^*=\sigma_i$, and,
following the convention just stated, $\sigma_i(P)$ means $\varphi_i(P)$,
so that $v_{\sigma_i(P)}(f)=v_P(\sigma_i(f))$ for every place $P$ and every
$f\in\F_5(X)^\times$. Since $\varphi_i$ and $\varphi_i^{-1}$ generate the same
group, this choice does not affect the orbits used below.

Throughout, $G=\langle\sigma_4\rangle\subseteq\Aut_{\F_5}(X)$ denotes the
cyclic group of order $12$ generated by the automorphism $\sigma_4$ given
explicitly in Appendix~\ref{app:aut} (Proposition~\ref{prop:autgroup}). All
results of this paper use only $G$. Magma's \texttt{Automorphisms} reports
that $G$ is the whole group $\Aut_{\F_5}(X)$ (Remark~\ref{rem:fullaut}), but
we do not use this anywhere: the reductions below remain valid for any group
of automorphisms, and a larger group would only allow further reduction.

\begin{lemma}\label{lem:aut}
Let $\sigma\in\Aut_{\F_5}(X)$ and $E\ge0$. Then $h\mapsto h\circ\sigma^{-1}$ is a
degree-preserving bijection $\LL(S+E)\to\LL(S+\sigma(E))$.
\end{lemma}

\begin{proof}
$\sigma$ permutes $X(\F_5)$, so $\sigma(S)=S$, and
$\divv(h\circ\sigma^{-1})=\sigma(\divv h)$.
\end{proof}

\section{The first search}\label{sec:tech}

The mathematical reductions that bring the search
down to about three million Riemann--Roch spaces are not new: they are all
taken from van Hoeij's code~\cite{vHX142} for $X_1(42)$, which accompanies
\cite{DvH} and proves $\gon_{\F_5}X_1(42)\ge 9$. Our contribution is the
organisation of the computation, above all its parallelization, which is what
makes the searches for $d=10$ and $d=11$ feasible. As noted earlier, we cannot guarantee that this search was exhaustive, nonetheless we record the method because the orbit-reduced search is a further refinement of this technique. 

\subsection{Symmetry reduction on rational points}\label{subsec:ratsym}
Let $O=\{Q_1,Q_2\}\subset X(\F_5)$, where $Q_1$ and $Q_2$ are the unique places
at which the function $w=x/y$ (the generator of $\F_5(X)$ over $\F_5(t)$ used
by Magma, see Appendix~\ref{app:aut}) has a zero of order $3$ and
a pole of order $7$, respectively. This is the choice made in~\cite{vHX142},
where $O$ is stated to contain one member of each diamond orbit on $X(\F_5)$.
The same reduction is used there to prove $\gon_{\F_5}X_1(42)\ge 9$.

\begin{proposition}\label{prop:orbits}
Every $G$-orbit in $X(\F_5)$ contains $Q_1$ or $Q_2$.
\end{proposition}

\begin{proof}
For the diamond group this is asserted in~\cite{vHX142}; since no proof is
given there, we verified it independently for $G$. For each
$\sigma=\varphi^*\in G$ and each rational place $P$ we identified the image
$\varphi(P)$ as the unique rational place $Q$ with
$v_Q(g)=v_P(\sigma(g))$ for a fixed list of test functions $g$. The resulting
permutation group on the twelve rational places is cyclic of order~$12$ and
acts simply transitively (Proposition~\ref{prop:autgroup}). In
particular there is a single orbit, and it contains $Q_1$ and $Q_2$.
\end{proof}

\begin{proposition}[Symmetry reduction]\label{prop:sym}
Let $B\ge0$ and let $\mathcal{E}$ be the set of effective divisors $E$ of degree
$B$ such that either $\supp E$ contains no rational point or
$\supp E\cap O\neq\varnothing$. Every effective divisor of degree $B$ is
$G$-equivalent to an element of $\mathcal{E}$.
\end{proposition}

\begin{proof}
If $\supp E$ contains a rational point $Q$, choose $\sigma\in G$ with
$\sigma(Q)\in O$ (Proposition~\ref{prop:orbits}). Then $\sigma(E)\in\mathcal{E}$.
\end{proof}

\begin{theorem}\label{thm:criterion}
Let $d\ge2$. Then $\gon_{\F_5}(X)\ge d+1$ if and only if, for every
$E\in\mathcal{E}$ of degree $d-2$, every nonconstant function in $\LL(S+E)$
has degree $>d$.
\end{theorem}

\begin{proof}
Combine Corollary~\ref{cor:criterion}, Proposition~\ref{prop:sym} and
Lemma~\ref{lem:aut}.
\end{proof}

\subsection{Enumeration}\label{subsec:enum}
An effective divisor of degree $B$ is a sum of places whose degrees form a
partition $\lambda\vdash B$, its \emph{type}. If $\lambda$ has $m_v$ parts equal
to $v$, the divisors of type $\lambda$ correspond to tuples of multisets of
$m_v$ places of degree $v$. To enumerate $\mathcal{E}$ we use all of them when
$m_1=0$. When $m_1\ge1$, we replace the multiset of rational places by one
element of $O$ together with a multiset of $m_1-1$ rational places. This
produces every element of $\mathcal{E}$, some more than once. The number of
enumeration instances for $\lambda$, that is, of Riemann--Roch spaces
computed, is therefore
\[
  N(\lambda)=\begin{cases}
  \prod_v\binom{P_v+m_v-1}{m_v}, & m_1=0,\\[1ex]
  2\binom{P_1+m_1-2}{m_1-1}\prod_{v\ge2}\binom{P_v+m_v-1}{m_v}, & m_1\ge1.
  \end{cases}
\]
Tables~\ref{tab:d10} and~\ref{tab:d11} list these numbers for $B=8$ and $B=9$.
These are counts of enumeration instances; since some divisors occur more
than once, the numbers of distinct divisors are slightly smaller. Without
the reduction of Section~\ref{subsec:ratsym} the totals would be the numbers
of all effective divisors, $1\,276\,671$ and $6\,370\,940$.

\begin{table}[ht]
\caption{Enumeration instances (Riemann--Roch spaces computed) for $d=10$ ($B=8$). Total: $600\,465$.}\label{tab:d10}
\small
\begin{tabular}{@{}lr@{\qquad}lr@{}}
\toprule
$\lambda$ & $N(\lambda)$ & $\lambda$ & $N(\lambda)$\\
\midrule
$[8]$ & 48\,708 & $[4\,1^{4}]$ & 133\,224\\
$[7\,1]$ & 22\,320 & $[3^{2}\,2]$ & 216\\
$[6\,2]$ & 13\,728 & $[3^{2}\,1^{2}]$ & 864\\
$[6\,1^{2}]$ & 54\,912 & $[3\,2^{2}\,1]$ & 336\\
$[5\,3]$ & 5\,184 & $[3\,2\,1^{3}]$ & 7\,488\\
$[5\,2\,1]$ & 7\,776 & $[3\,1^{5}]$ & 21\,840\\
$[5\,1^{3}]$ & 101\,088 & $[2^{4}]$ & 126\\
$[4^{2}]$ & 16\,836 & $[2^{3}\,1^{2}]$ & 1\,344\\
$[4\,3\,1]$ & 2\,928 & $[2^{2}\,1^{4}]$ & 15\,288\\
$[4\,2^{2}]$ & 3\,843 & $[2\,1^{6}]$ & 52\,416\\
$[4\,2\,1^{2}]$ & 26\,352 & $[1^{8}]$ & 63\,648\\
\bottomrule
\end{tabular}
\end{table}

\begin{table}[ht]
\caption{Enumeration instances (Riemann--Roch spaces computed) for $d=11$ ($B=9$). Total: $2\,994\,534$.}\label{tab:d11}
\small
\begin{tabular}{@{}lr@{\qquad}lr@{}}
\toprule
$\lambda$ & $N(\lambda)$ & $\lambda$ & $N(\lambda)$\\
\midrule
$[9]$ & 214\,384 & $[4\,2^{2}\,1]$ & 7\,686\\
$[8\,1]$ & 97\,416 & $[4\,2\,1^{3}]$ & 171\,288\\
$[7\,2]$ & 66\,960 & $[4\,1^{5}]$ & 499\,590\\
$[7\,1^{2}]$ & 267\,840 & $[3^{3}]$ & 120\\
$[6\,3]$ & 18\,304 & $[3^{2}\,2\,1]$ & 432\\
$[6\,2\,1]$ & 27\,456 & $[3^{2}\,1^{3}]$ & 5\,616\\
$[6\,1^{3}]$ & 356\,928 & $[3\,2^{3}]$ & 448\\
$[5\,4]$ & 118\,584 & $[3\,2^{2}\,1^{2}]$ & 4\,032\\
$[5\,3\,1]$ & 10\,368 & $[3\,2\,1^{4}]$ & 34\,944\\
$[5\,2^{2}]$ & 13\,608 & $[3\,1^{6}]$ & 69\,888\\
$[5\,2\,1^{2}]$ & 93\,312 & $[2^{4}\,1]$ & 252\\
$[5\,1^{4}]$ & 471\,744 & $[2^{3}\,1^{3}]$ & 8\,736\\
$[4^{2}\,1]$ & 33\,672 & $[2^{2}\,1^{5}]$ & 57\,330\\
$[4\,3\,2]$ & 8\,784 & $[2\,1^{7}]$ & 148\,512\\
$[4\,3\,1^{2}]$ & 35\,136 & $[1^{9}]$ & 151\,164\\
\bottomrule
\end{tabular}
\end{table}

\subsection{Testing a divisor}\label{subsec:test}
For each $E$ we call the Magma function \texttt{RiemannRochSpace} on
$S+E$, which implements Hess's algorithm~\cite{Hess}. It returns an
$\F_5$-vector space $V$ of dimension $\ell(S+E)$ together with an
$\F_5$-linear isomorphism $\mu\colon V\to\LL(S+E)\subset\F_5(X)$. Since
$\LL(S+E)$ contains the constants, $\ell(S+E)\ge1$; if $\ell(S+E)=1$, the
space consists of constants only and nothing further is done. Otherwise the
program runs through all $5^{\ell(S+E)}$ vectors $v\in V$, discards the five
vectors with $\mu(v)\in\F_5$, and computes $\deg\mu(v)$, the degree of the pole
divisor of $\mu(v)$, for each of the remaining $5^{\ell(S+E)}-5$ functions.
Thus every nonconstant element of $\LL(S+E)$ is tested, not only the basis
elements; this matters because a linear combination of basis elements can
have fewer poles than the basis elements themselves. Testing one element of
each class $\{ah+b: a\in\F_5^\times,\ b\in\F_5\}$ would suffice, since
$\deg(ah+b)=\deg h$, but we did not use this shortcut. A hit is recorded
whenever $\deg\mu(v)\le d$. Functions of degree $\le 12$ are also
written to a log with their pole divisors.

In both searches the log contains functions of degree $12$, and no function
of degree $\le 11$. This is consistent with \eqref{eq:upper}: since
$\deg(S+E)=n+d-2\ge 20$, the spaces $\LL(S+E)$ can contain functions of degree
$12$, and the search finds them. Van Hoeij's computation~\cite{vHX142} also
found a function of degree $12$ over $\F_5$, in a space of the form
$\LL(S+Q+2Q'+3Q'')$ with $Q,Q'$ rational and $Q''\in O$. The logged
functions give an independent witness of $\gon_{\F_5}X_1(42)\le 12$.

\subsection{Parallel computation}\label{subsec:20cpu}
Van Hoeij's code is a single sequential computation, which suffices for the
search up to degree $8$ carried out in~\cite{vHX142}. For $d=10$ and $d=11$ the number of divisors and the cost of the
places of degree $8$ and $9$ make a sequential run impractical. We
reorganised the search so that it splits into independent jobs. Both searches were run on a single workstation with an Intel Xeon W-2265
processor ($3.5$\,GHz, $12$ physical cores, $24$ hardware threads) and
$64$\,GB of memory, as $20$ independent Magma processes (Magma V2.29-10), each occupying one of
the $24$ logical CPUs.

For $d=10$, job $r\in\{0,\dots,19\}$ tests exactly the divisors whose index
$i$ satisfies $i\equiv r+1\pmod{20}$, over all partitions of $8$. The residue
classes partition each index set, so the union of the jobs is the full
enumeration, and the round-robin assignment interleaves cheap and expensive
divisors, so the jobs finish at about the same time.

For $d=11$ the same scheme exhausted the memory: every job had to hold the
$214\,384$ places of degree $9$ and the $48\,708$ places of degree $8$. We
therefore split the $30$ partitions of $9$ into three groups,
\[
  \{[9]\},\qquad \{[8,1]\},\qquad \text{the remaining } 28 \text{ partitions},
\]
with $1$, $1$ and $18$ jobs respectively, and applied
the round-robin rule within each group. Each job computes only the places of
the degrees occurring in its own group, so only one job ever stores the
places of degree $9$. The partitions $[9]$ and $[8,1]$ account for only about
$10\%$ of the divisors ($311\,800$ of $2\,994\,534$), but most of the memory.

We used only runs in which all $20$ jobs reported
\texttt{COMPLETED} and all jobs of a group reported identical place counts and
partition sizes. These checks do not, however, certify that the division of
the work is a partition: the index $i$ of a divisor depends on the order in
which Magma returns the places of each degree, and this order is not
reproducible between Magma sessions (see Section~\ref{subsec:split}). If two
jobs received the places in different orders, some divisors could have been
tested twice and others not at all. We therefore do not use the searches of
this section in the proof of the Theorem; they serve as supporting evidence,
consistent with the orbit-reduced search of Section~\ref{sec:orbitsearch},
whose division of the work does not depend on this order.

For $d=11$ the $20$ jobs used $11\,973\,454$ seconds of CPU time in total
(about $3\,326$ hours, or $139$ days), while the longest job used
$629\,624$ seconds (about $7.3$ days); these figures count the search itself,
after the places have been computed. The parallel run therefore completed in
about a week instead of the more than four months a sequential run would have
needed, a speed-up of about $19$. For $d=10$ the corresponding figures are
$2\,409\,708$ seconds in total (about $669$ hours, or $28$ days) and
$125\,522$ seconds (about $1.5$ days) for the longest job, again a speed-up
of about $19$.

\section{An orbit-reduced search}\label{sec:orbitsearch}

In this section we describe a separately implemented exhaustive
search for $d=11$. It uses a
different enumeration of divisors from the search of Section~\ref{sec:tech},
but shares with it the model of $X$ over $\F_5$, the search criterion of
Section~\ref{sec:criterion}, Magma's Riemann--Roch machinery and the
automorphism $\sigma_4$. It
uses the action of $G$ on the places of degree $\le5$, not
only on the rational places. Each divisor $E$ of degree $9$ is split as
$E=E_{\le5}+E_{\ge6}$ into its parts supported on places of degree $\le5$ and
$\ge6$; the search keeps one representative of each $G$-orbit of parts
$E_{\le5}$, combined with every possible part $E_{\ge6}$.  What is
new here is the specific use of the action on places of degree $2,\dots,5$ and
its implementation. For comparison, in van Hoeij's code~\cite{vHX142} for
$N=42$ the automorphisms enter only through the set $O$ of
Section~\ref{subsec:ratsym}, that is, only through their action on the
rational places, and places of higher degree are enumerated without any
reduction. We are not aware of an earlier implementation for $X_1(N)$ that
uses the action on places of higher degree, and the algorithm of~\cite{DT}
restricts the divisors by a counting argument on the rational points
instead.

\subsection{The group $G$ on places of small degree}
Write $G=\langle\sigma\rangle$, $\sigma=\sigma_4$, for the cyclic group of
automorphisms of order $12$ of Section~\ref{sec:sym}
(Proposition~\ref{prop:autgroup}). By Section~\ref{sec:prelim},
$G$ acts on the places of each degree $v$, and since $G$ is cyclic its orbits
are the cycles of the permutation induced by $\sigma$.

For every element of $G$ we also computed its pull-back, in particular $\sigma^*$. The image $\sigma(P)$ of a place $P$ is determined by
$v_{\sigma(P)}(f)=v_P(\sigma^*f)$ for all $f\in\F_5(X)^\times$, so it can be
computed exactly from the pull-back data (see the data file \texttt{X1\_42\_aut\_F5.txt}).
Equivalently, if $g\in\F_5(X)$ is regular at $\sigma(P)$, then $\sigma^*g$ is
regular at $P$ and
\[
  g\bigl(\sigma(P)\bigr)=(\sigma^*g)(P)\qquad\text{in }\kappa(P)\cong\kappa\bigl(\sigma(P)\bigr),
\]
so the values of the coordinate functions at $\sigma(P)$ are the values of
their pull-backs at $P$.

\begin{example}
The following toy case on $\PP^1_{\F_5}$, with function field $\F_5(t)$,
shows the computation. Let $\sigma\colon\PP^1\to\PP^1$ be $a\mapsto a+1$, so
that $\sigma^*t=t\circ\sigma=t+1$.
\begin{enumerate}[label=(\alph*),leftmargin=*]
  \item Let $P$ be the rational place $t=2$, that is, the zero of $t-2$.
    Then $t(\sigma(P))=(\sigma^*t)(P)=2+1=3$, so $\sigma(P)$ is the place
    $t=3$. In terms of valuations: $v_P\bigl(\sigma^*(t-3)\bigr)=v_P(t-2)=1$,
    so $v_{\sigma(P)}(t-3)=1$.
  \item Let $P$ be the place of degree $2$ given by the irreducible
    polynomial $t^2+2$ (note that $-2=3$ is not a square modulo $5$), and let
    $\alpha\in\kappa(P)\cong\F_{25}$ be the residue of $t$, so
    $\alpha^2+2=0$. Then $t(\sigma(P))=(\sigma^*t)(P)=\alpha+1$, whose minimal
    polynomial over $\F_5$ is
    \[
      (u-1)^2+2=u^2+3u+3 .
    \]
    Hence $\sigma(P)$ is the place of degree $2$ given by $t^2+3t+3$.
    Indeed, $\sigma^*(t^2+3t+3)=(t+1)^2+3(t+1)+3=t^2+5t+7=t^2+2$ in $\F_5[t]$,
    so $v_P\bigl(\sigma^*(t^2+3t+3)\bigr)=1$ and
    $v_{\sigma(P)}(t^2+3t+3)=1$.
\end{enumerate}
\end{example}

On $X$, with $\F_5(X)=\F_5(t)[w]/(\Phi(t,w))$ as in Appendix~\ref{app:aut}
(so $t=1/y$ and $w=x/y$ in the coordinates of the model), we determine $\sigma(P)$ for a place $P$ of degree $v$ in three steps.
First, the minimal polynomial over $\F_5$ of $(\sigma^*t)(P)$ (or the fact
that $\sigma^*t$ has a pole at $P$) gives the place of $\F_5(t)$ below
$\sigma(P)$, and the candidates are the places of degree $v$ above it. Second, if there are several candidates and both $\sigma^*t$ and $\sigma^*w$
are regular at $P$, then $t$ and $w$ are regular at $\sigma(P)$, hence so are
$t+w$ and $t+2w$, and we compare the minimal polynomials of the values of
$t$, $w$, $t+w$ and $t+2w$ at $\sigma(P)$, computed from the pull-backs
$\sigma^*t$, $\sigma^*w$, $\sigma^*(t+w)$, $\sigma^*(t+2w)$ at $P$, with those
of the candidates. A candidate at which $t$ or $w$ has a pole cannot match, therefore, if $\sigma^*t$ or $\sigma^*w$ has a pole at $P$, this step is skipped. Third, if several candidates still remain, we take for each
candidate $Q$ a local uniformizer $u_Q$ at $Q$; then $u_Q$ vanishes at
$\sigma(P)$ if and only if $v_P(\sigma^*u_Q)>0$, and this vanishing pattern
determines $\sigma(P)$. In our computation these steps determined $\sigma(P)$
for every place of degree $\le5$, and the images form a permutation of each
set of places.

In this way we computed the
permutation of $\sigma$ on the places of each degree $v\le5$; the resulting
orbits are listed in Appendix~\ref{app:orbits}.

\subsection{The reduction}
Every effective divisor $E$ of degree $9$ decomposes uniquely as
\[
  E=E_{\le5}+E_{\ge6},
\]
where $E_{\le5}$ is supported on places of degree $\le5$ and $E_{\ge6}$ on
places of degree $\ge6$. Since $2\cdot6>9$, the part $E_{\ge6}$ is either $0$
or a single place. The group $G$ acts on the set of effective divisors
supported on places of degree $\le5$. Choose one representative in each
$G$-orbit of this set, and let $\mathcal{R}$ be the set of chosen
representatives.

\begin{theorem}\label{thm:orbitcriterion}
Let $\mathcal{E}^{\mathrm{orb}}$ be the set of effective divisors $E$ of
degree $9$ with $E_{\le5}\in\mathcal{R}$. Then $\gon_{\F_5}(X)\ge12$ if and only if, for every
$E\in\mathcal{E}^{\mathrm{orb}}$, every nonconstant function in $\LL(S+E)$
has degree $>11$.
\end{theorem}

\begin{proof}
Let $E$ be any effective divisor of degree $9$. The orbit of $E_{\le5}$
contains exactly one element of $\mathcal{R}$, say $\sigma^k(E_{\le5})$.
Since $\sigma$ preserves degrees of places,
$\sigma^k(E)_{\le5}=\sigma^k(E_{\le5})\in\mathcal{R}$, so
$\sigma^k(E)\in\mathcal{E}^{\mathrm{orb}}$.
By Lemma~\ref{lem:aut}, $\LL(S+E)$ and $\LL(S+\sigma^k(E))$ contain functions
of the same degrees. The statement now follows from
Corollary~\ref{cor:criterion} with $d=11$.
\end{proof}

The proof uses neither the action of $\sigma$ on places of degree $\ge6$ nor
the set $O$. For divisors with $E_{\ge6}=0$, $\mathcal{E}^{\mathrm{orb}}$
contains exactly one divisor from each $G$-orbit of divisors. In general,
it contains the representative $E_{\le5}\in\mathcal{R}$ of each $G$-orbit of
low-degree parts, combined with every place of degree $\ge6$; it need not
contain exactly one divisor from each $G$-orbit of divisors $E$. In the computation $\mathcal{R}$ is chosen by the following rule. The places of degree
$\le5$ carry labels that do not depend on the Magma session
(Section~\ref{subsec:split}), and the labels of each degree are totally
ordered. We write a part $E_{\le5}$ as the tuple, over the degrees occurring
in it in decreasing order, of the sorted lists of labels of its places, each
place listed as often as its multiplicity. A part $E_{\le5}$ belongs to
$\mathcal{R}$ if and only if this tuple is lexicographically smallest among
the tuples of $\sigma^k(E_{\le5})$, $k=0,\dots,11$. Since the permutations of
$\sigma$ on the labels are stored, $\mathcal{R}$ is computed without any
further computation on the curve.

For each partition $\lambda\vdash9$, the number of admissible $E_{\le5}$ was
computed by enumeration. As a check, we compared it with the number of
$G$-orbits obtained from Burnside's lemma using only the cycle types of
$\sigma$ (Appendix~\ref{app:orbits}).
The two agree for every $\lambda$. Table~\ref{tab:orb} lists the numbers. In
total $743\,231$ divisors are tested, compared with $2\,994\,534$ in
Table~\ref{tab:d11} and $6\,370\,940$ without any reduction.

\begin{table}[ht]
\caption{The orbit-reduced search for $d=11$: the number $n_{\le5}$ of
admissible parts $E_{\le5}$, and the number of divisors tested. Total:
$743\,231$.}\label{tab:orb}
\small
\begin{tabular}{@{}lrr@{\qquad}lrr@{}}
\toprule
$\lambda$ & $n_{\le5}$ & divisors & $\lambda$ & $n_{\le5}$ & divisors\\
\midrule
$[9]$ & 1 & 214\,384 & $[4\,3\,2]$ & 732 & 732\\
$[8\,1]$ & 1 & 48\,708 & $[4\,3\,1^2]$ & 9\,516 & 9\,516\\
$[7\,2]$ & 1 & 11\,160 & $[4\,2^2\,1]$ & 3\,843 & 3\,843\\
$[7\,1^2]$ & 7 & 78\,120 & $[4\,2\,1^3]$ & 33\,306 & 33\,306\\
$[6\,3]$ & 2 & 4\,576 & $[4\,1^5]$ & 66\,612 & 66\,612\\
$[6\,2\,1]$ & 6 & 13\,728 & $[3^3]$ & 30 & 30\\
$[6\,1^3]$ & 31 & 70\,928 & $[3^2\,2\,1]$ & 216 & 216\\
$[5\,4]$ & 9\,882 & 9\,882 & $[3^2\,1^3]$ & 1\,116 & 1\,116\\
$[5\,3\,1]$ & 5\,184 & 5\,184 & $[3\,2^3]$ & 40 & 40\\
$[5\,2^2]$ & 1\,134 & 1\,134 & $[3\,2^2\,1^2]$ & 1\,092 & 1\,092\\
$[5\,2\,1^2]$ & 25\,272 & 25\,272 & $[3\,2\,1^4]$ & 5\,460 & 5\,460\\
$[5\,1^4]$ & 73\,710 & 73\,710 & $[3\,1^6]$ & 8\,264 & 8\,264\\
$[4^2\,1]$ & 16\,836 & 16\,836 & $[2^4\,1]$ & 126 & 126\\
 & & & $[2^3\,1^3]$ & 1\,700 & 1\,700\\
 & & & $[2^2\,1^5]$ & 7\,644 & 7\,644\\
 & & & $[2\,1^7]$ & 15\,912 & 15\,912\\
 & & & $[1^9]$ & 14\,000 & 14\,000\\
\bottomrule
\end{tabular}
\end{table}

The partition $[9]$ is not reduced at all, since it has no part of degree
$\le5$. Reducing it would require the orbits of $G$ on the $214\,384$ places of
degree $9$. For $[8,1]$ the reduction is by the factor $12$ from the rational
place, while the old method gives a factor of $6$.

\subsection{Splitting the work}\label{subsec:split}
The divisors are divided into five groups according to the largest part of
$\lambda$: $[9]$, $[8,1]$, the partitions with largest part $7$, those with
largest part $6$, and all others (largest part $\le5$). The admissible parts
$E_{\le5}$ are computed once, before the search.
\begin{itemize}[leftmargin=*]
  \item In the last group a divisor is just its part $E_{\le5}$, and the
    divisors are dealt out to the jobs in turn.
  \item In a group with high degree $v\in\{6,7,8,9\}$, each of the jobs
    receives all admissible $E_{\le5}$ of the group, and the places of degree
    $v$ are divided among the jobs so that each place is handled by exactly
    one job. A job tests $S+E_{\le5}+Q$ for all its $E_{\le5}$ and all
    places $Q$ in its share.
\end{itemize}
Magma does not return the places of a given degree $v$ in the same order in
different sessions (\texttt{Places(AFF,}$v$\texttt{)}), so neither the
stored parts $E_{\le5}$ nor the shares of the high-degree places can be
described by positions in such a list. We therefore identify places by data
that depend only on the place itself. Every place $Q$ of degree $\ge2$ is
assigned its \emph{bucket} $b(Q)$, the minimal polynomial over $\F_5$ of the
value of $t$ at $Q$ (with a separate bucket for the poles of $t$).

\emph{Labels of the places of degree $\le5$.} A place of degree $2\le v\le5$
is described by its bucket together with the minimal polynomials over $\F_5$
of the values of $w$, $t+w$ and $t+2w$ at it. We checked that these data are
distinct for distinct places of each degree $2\le v\le5$; the label of a
place is the rank of its data in a fixed sorted order. The twelve rational
places do not all have distinct data of this kind, and are labelled instead
by their \emph{signatures}: the valuations, and the values where the
valuation is $0$, of a fixed list of functions ($t$, $w$, a local
uniformizer at each rational place and a basis of $\LL(3S)$), which is
stored in the data repository as explicit elements of $\F_5(X)$; we checked
that the twelve signatures are distinct. Each job recomputes the places of
degree $\le5$ that it needs, recomputes their labels, and stops with an
error unless every label it needs is found exactly once.

\emph{Shares of the high-degree places.} Let a group with high degree $v$ have
$s$ jobs, numbered $0,\dots,s-1$. Job $r$ computes all places of degree $v$,
checks that there are $P_v$ of them, and keeps exactly the places $Q$ with
$H(b(Q))\equiv r\pmod s$, where $H$ is the following fixed function of the
bucket. A monic polynomial $g=\sum_{i=0}^{m}g_iu^i\in\F_5[u]$ of degree
$m$, with $g_i\in\{0,\dots,4\}$ (so $g_m=1$), is encoded as the integer
$\operatorname{enc}(g)=5^{m+1}+\sum_{i=0}^{m}g_i5^i$, and the bucket of the
poles of $t$ as $0$; then
\[
  H(b)=\bigl((\operatorname{enc}(b)\bmod N_0)\cdot 1\,000\,003+12\,345\bigr)\bmod N_0,
  \qquad N_0=10^9+7 .
\]
Any fixed function of the bucket would do; this one only serves to spread
the buckets evenly over the jobs. Since $b(Q)$ depends only on $Q$, every place of degree $v$ is
kept by exactly one job of the group, whatever order Magma returns the
places in. As a further check, each job reports a fingerprint (a sum of
hash values modulo a large prime) of the buckets of all its places of degree
$v$ and of the buckets of the places it kept; we checked that all jobs of a
group report the same fingerprint for all places, that the kept
fingerprints add up to it, and that the numbers of kept places add up to
$P_v$. In the last group the divisors are dealt out to the jobs in turn when
the job files are written, once and in a single process, so this division
does not depend on Magma's order either.

The run is accepted only if every job reports \texttt{COMPLETED}, all the
checks above are satisfied, and the numbers of divisors tested add up to
$743\,231$. The scripts and checks are documented in the data repository.

The search was run with Magma V2.29-10 on the same workstation as the search
of Section~\ref{sec:tech}: an Intel Xeon W-2265 with
$12$ physical cores, $24$ hardware threads and $64$\,GB of memory, again as
$20$ independent Magma processes: $2$ jobs for $[9]$, $2$ for $[8,1]$, $3$ each
for the partitions with largest part $7$ and $6$, and $10$ for the remaining
partitions. The number of jobs for $[9]$ was limited by memory: each of them
holds all $214\,384$ places of degree $9$, which took about $4\,900$ seconds
and about $19$\,GB of memory to compute. Each space $\LL(S+E)$ was tested exactly as described in
Section~\ref{subsec:test}. All $20$ jobs reported
\texttt{COMPLETED}, all the checks above were satisfied, and no space
$\LL(S+E)$ contained a nonconstant function of degree $\le11$; the smallest
degree found was $12$. The $20$ jobs used $2\,077\,413$ seconds of CPU time in
total (about $577$ hours, or $24$ days), and the longest job used
$164\,850$ seconds (about $1.9$ days); as in Section~\ref{sec:tech}, these
figures count the search itself, after the places have been computed. Compared
with the search of Section~\ref{sec:tech} ($11\,973\,454$ and $629\,624$
seconds), the total CPU time is smaller by a factor of about $5.8$ and the
longest job by a factor of about $3.8$.

\section{Proof of the theorem}\label{sec:proof}

\begin{proof}[Proof of the Theorem]
The orbit-reduced search of Section~\ref{sec:orbitsearch} ran over all
$E\in\mathcal{E}^{\mathrm{orb}}$ (Table~\ref{tab:orb}), and no
$\LL(S+E)$ contained a nonconstant function of degree $\le 11$.
By Theorem~\ref{thm:orbitcriterion}, $\gon_{\F_5}X_1(42)\ge 12$.

The searches of Section~\ref{sec:tech} are consistent with this: with
$d=11$ they found no nonconstant function of degree $\le11$ in any of the
spaces they tested. As explained in Section~\ref{subsec:20cpu}, we cannot
certify that they tested every $E\in\mathcal{E}$, so we use them only as
supporting evidence. They share with the orbit-reduced search the model of
$X$ over $\F_5$ (Lemma~\ref{lem:model}), Corollary~\ref{cor:criterion}, the
automorphism $\sigma_4$ and Magma's computation of places and Riemann--Roch
spaces, so their agreement is a check on the two enumerations and their
implementations, not on these shared ingredients.

Together with \eqref{eq:red} and \eqref{eq:upper} this gives
\[
  12\le\gon_{\F_5}X_1(42)\le\gon_\Q X_1(42)\le12 ,
\]
so all three quantities are equal.
\end{proof}

\section*{Data availability}
The Magma scripts and the log files of the computation are available at
\url{https://github.com/soura-code/X1_42_gonality}. They contain, for every
job, the place counts, the partition sizes, the number of divisors tested
and the final status. The repository also contains Sutherland's model
\texttt{FFFc42.txt}, the script \texttt{X1\_42\_moduli\_check.m} that carries
out the verifications in the proof of Lemma~\ref{lem:model}, the data file
\texttt{X1\_42\_aut\_F5.txt} with the twelve automorphisms of
$X_1(42)_{\F_5}$, and the file
\texttt{X1\_42\_\allowbreak automorphisms.tex}, which lists them explicitly together with
their orders, their expressions as powers of $\sigma_4$ and their permutations
of the rational places (Appendix~\ref{app:aut}).

\section*{Acknowledgements}
I am deeply grateful to my advisor, Kyu-Hwan Lee, for his guidance,
encouragement and constant support throughout this project. I thank the Connecticut Number Theory Research Group, and in
particular \'Alvaro Lozano-Robledo, for providing access to Anthropic's
Claude Pro, which has been instrumental. (see the \hyperref[sec:ai]{Declaration of generative AI use})

\section*{Declaration of generative AI use}\label{sec:ai}

During the preparation of this manuscript, the author used Anthropic's
Claude to explore mathematical ideas and proof strategies, to develop and
debug the Magma and Python code used in the computations, to search the
literature, and to assist with editorial and formatting checks. The author
reviewed and verified all AI-generated suggestions, code and computational
results.

\appendix
\section{The automorphism group of $X_1(42)_{\F_5}$}\label{app:aut}

\subsection{Conventions}
Magma represents the function field of $X=X_1(42)_{\F_5}$ (built from the
\texttt{FFFc42} model as in Section~\ref{sec:prelim}) as a simple extension
\[
  \F_5(X)=\F_5(t)[w]/\bigl(\Phi(t,w)\bigr),\qquad [\F_5(X):\F_5(t)]=16,
\]
where $t$ generates the rational subfield (\texttt{BaseField(AFF).1}) and $w$
is the generator \texttt{AFF.1}. Magma builds this field from the projective closure of $\bar f(x,y)=0$ in
the affine chart where $y\ne0$, so that, in terms of the coordinates $x,y$ of
the plane model of Section~\ref{sec:prelim},
\[
  t=\frac1y,\qquad w=\frac xy,\qquad
  \Phi(t,w)=t^{23}\,\bar f\Bigl(\frac wt,\frac1t\Bigr),
\]
where $23$ is the total degree of $\bar f$. The polynomial $\Phi$ is
\texttt{DefiningPolynomial(AFF)}; it has degree $16=\deg_x\bar f$ in $w$,
with leading coefficient $t^3$, and it is included in the data repository.
An $\F_5$-automorphism $\sigma$ of $\F_5(X)$ is determined by $\sigma(t)$ and
$\sigma(w)$, and every element of $\F_5(X)$ is written uniquely as
$\sum_{j=0}^{15}c_j(t)\,w^{j}$ with $c_j\in\F_5(t)$.

The call \texttt{Automorphisms(AFF)} in Magma returns $12$ maps
$\sigma_1=\mathrm{id},\sigma_2,\dots,\sigma_{12}$ (computing them took
$5814$ seconds). We checked that each $\sigma_i$ respects the equation
$\Phi(t,w)=0$, that is, $\Phi(\sigma_i(t),\sigma_i(w))=0$,
so that it is an $\F_5$-endomorphism of $\F_5(X)$, that $\sigma_1$ is the
identity, and that the set $\{\sigma_1,\dots,\sigma_{12}\}$ is closed under
composition.
For each $\sigma_i=\varphi_i^*$ we determined the permutation of $X(\F_5)$
induced by $\varphi_i$ as in the proof of Proposition~\ref{prop:orbits}. The order of each $\sigma_i$, its expression as a power of $\sigma_4$ and its
permutation of the rational places are listed in the file
\texttt{X1\_42\_\allowbreak automorphisms.tex} in the data repository; there
the rational places are labelled by their signatures
(Section~\ref{subsec:split}), which do not depend on the Magma session.

\begin{proposition}\label{prop:autgroup}
\begin{enumerate}
  \item The maps $\sigma_1,\dots,\sigma_{12}$ are automorphisms of $X$, and
        they form a cyclic group $G=\langle\sigma_4\rangle$ of order $12$.
  \item $G$ acts simply transitively on the twelve rational places, and
        $\sigma_4^5(Q_2)=Q_1$.
  \item The subgroup $C=\langle\sigma_4^2\rangle$ of order $6$ has two orbits
        on $X(\F_5)$, each of size $6$, and $Q_1$ and $Q_2$ lie in different
        orbits.
\end{enumerate}
\end{proposition}

\begin{proof}
Let $M=\{\sigma_1,\dots,\sigma_{12}\}$. Each $\sigma\in M$ is an injective
$\F_5$-endomorphism of $\F_5(X)$. Since $M$ is finite and closed under
composition, $\sigma^m=\sigma^n$ for some $m>n\ge1$, and injectivity gives
$\sigma^{m-n}=\mathrm{id}$; so $\sigma$ is invertible with inverse
$\sigma^{m-n-1}\in M$. Hence $M$ is a group of automorphisms of $\F_5(X)$,
that is, of $X$, of order $12$. We computed the permutation of $\sigma_4$ on
the twelve rational places, identified by their signatures
(Section~\ref{subsec:split}): it is a single cycle of length $12$, and
$\sigma_4^5(Q_2)=Q_1$. Hence $\sigma_4$ has order at least $12$, so
$M=\langle\sigma_4\rangle$; the map from $M$ to the symmetric group on
$X(\F_5)$ is injective, and its image, generated by a $12$-cycle, acts simply
transitively. The orbits of $C=\langle\sigma_4^2\rangle$ consist of the places
$\sigma_4^{k}(Q_2)$ with $k$ even and with $k$ odd, respectively; since $5$
is odd, $Q_1=\sigma_4^5(Q_2)$ and $Q_2$ lie in different orbits.
\end{proof}

\begin{remark}[The full automorphism group]\label{rem:fullaut}
The Magma function \texttt{Automorphisms}, called without any
restriction as \texttt{Automorphisms(AFF)}, is documented to return all automorphisms of the function
field. Its output consists of exactly the $12$ maps above, so according to
this computation $G=\Aut_{\F_5}(X)$; in particular $\Aut_{\F_5}(X)$ is cyclic
of order $12$. None of the results of this paper uses
$G=\Aut_{\F_5}(X)$: all of them are proved for the group $G$ of
Proposition~\ref{prop:autgroup}.
\end{remark}

\begin{remark}\label{rem:diamond-orbits}
Because the action of $G$ on $X(\F_5)$ is transitive, the set $O$
in Section~\ref{subsec:ratsym} could be replaced by $\{Q_1\}$ alone; we did not use
this further reduction.
\end{remark}

\subsection{Orbits on places of small degree}\label{app:orbits}
The permutations of $\sigma=\sigma_4$ on the places of degree $v\le5$ used in
Section~\ref{sec:orbitsearch} have the following orbits.
\begin{center}
\begin{tabular}{@{}cccl@{}}
\toprule
$v$ & $P_v$ & orbits & orbit sizes\\
\midrule
$1$ & $12$ & $1$ & $12$\\
$2$ & $6$ & $1$ & $6$\\
$3$ & $8$ & $2$ & $4,4$\\
$4$ & $183$ & $17$ & $3,\ 6^2,\ 12^{14}$\\
$5$ & $648$ & $54$ & $12^{54}$\\
\bottomrule
\end{tabular}
\end{center}
All orbit sizes divide $12$. In degree $4$ the group does not act freely: one
orbit consists of places fixed by $\sigma^3$, and two of places fixed by
$\sigma^6$.

\subsection{The generator $\sigma_4$}
Each $\sigma_i(t)$ and $\sigma_i(w)$ is an element
$\sum_{j=0}^{15}c_j(t)\,w^j$ of $\F_5(X)$ with $c_j\in\F_5(t)$. After
cancellation, every denominator that occurs is of the form
\[
  t^{a}\,(t+4)^{b}\,(t^2+4t+1)^{c}=t^{a}(t-1)^{b}(t^2-t+1)^{c}.
\]
We therefore write $c_j=N_j(t)/\bigl(t^{a_j}(t+4)^{b_j}(t^2+4t+1)^{e_j}\bigr)$
and record, for each $j$, the exponents $a_j,b_j,e_j$ and the coefficients of
the numerator $N_j\in\F_5[t]$, from the leading coefficient down to the
constant term; a row with entry $0$ means $c_j=0$. For example, a row
$j\;|\;2\;1\;0\;|\;2\;|\;3\,0\,1$ stands for
$c_j=(3t^2+1)/\bigl(t^2(t+4)\bigr)$.

By Proposition~\ref{prop:autgroup}, $G$ is cyclic of order $12$,
generated by $\sigma_4$, so the whole group is determined by $\sigma_4$: every
other element is a power $\sigma_4^k$, with $k$ given in Table~1 of the file \texttt{X1\_42\_\allowbreak automorphisms.tex}, which also contains the explicit formulas for all twelve automorphisms. We therefore list only the generator. The tables are generated mechanically from
the data file \texttt{X1\_42\_aut\_F5.txt} by the script
\texttt{aut\_to\_latex.py}, which checks the factorization of every
denominator exactly; they are an exact transcription of Magma's output.

\begin{center}\scriptsize
$\sigma_{4}(t)=\displaystyle\sum_{j=0}^{15}\frac{N_j(t)}{t^{a_j}(t+4)^{b_j}(t^2+4t+1)^{e_j}}\,w^j$
\\[4pt]
\begin{tabular}{@{}rrrrrl@{}}
\toprule
$j$ & $a_j$ & $b_j$ & $e_j$ & $\deg N_j$ & coefficients of $N_j$, from $t^{\deg N_j}$ down to $t^0$\\
\midrule
$0$ & $17$ & $10$ & $0$ & $25$ & $4\,1\,1\,0\,3\,0\,3\,3\,3\,2\,1\,1\,2\,4\,4\,1\,4\,1\,2\,3\,4\,4\,3\,4\,2\,4$\\
$1$ & $18$ & $11$ & $0$ & $28$ & $4\,0\,1\,2\,2\,1\,4\,2\,4\,3\,0\,1\,0\,4\,4\,0\,0\,3\,1\,4\,3\,2\,2\,0\,4\,4\,4\,3\,2$\\
$2$ & $19$ & $11$ & $0$ & $26$ & $4\,1\,0\,4\,3\,1\,2\,2\,2\,0\,0\,4\,4\,2\,3\,3\,1\,3\,4\,4\,4\,2\,3\,3\,2\,4\,3$\\
$3$ & $20$ & $10$ & $0$ & $27$ & $1\,4\,0\,0\,2\,1\,4\,4\,0\,2\,4\,4\,2\,4\,1\,0\,0\,0\,4\,4\,4\,1\,3\,4\,4\,4\,0\,2$\\
$4$ & $21$ & $12$ & $0$ & $28$ & $3\,1\,4\,3\,3\,1\,4\,4\,3\,2\,4\,3\,4\,2\,3\,4\,2\,0\,1\,3\,0\,2\,0\,1\,1\,1\,0\,3\,4$\\
$5$ & $27$ & $15$ & $0$ & $37$ & $4\,1\,0\,4\,0\,3\,3\,0\,3\,4\,4\,2\,1\,4\,4\,1\,0\,2\,2\,0\,4\,4\,4\,4\,0\,0\,4\,3\,3\,4\,1\,3\,1\,4\,3\,4\,2\,4$\\
$6$ & $28$ & $16$ & $0$ & $37$ & $1\,4\,3\,2\,2\,2\,2\,3\,0\,1\,3\,2\,1\,0\,2\,2\,0\,4\,3\,1\,2\,4\,1\,2\,2\,2\,0\,3\,4\,2\,2\,4\,1\,4\,4\,4\,3\,2$\\
$7$ & $29$ & $16$ & $0$ & $35$ & $3\,1\,2\,3\,1\,4\,3\,4\,0\,1\,3\,3\,1\,1\,2\,2\,4\,0\,4\,3\,2\,4\,0\,4\,4\,3\,3\,1\,4\,2\,2\,2\,3\,3\,3\,3$\\
$8$ & $30$ & $16$ & $0$ & $36$ & $3\,1\,2\,3\,2\,3\,3\,1\,4\,2\,1\,1\,1\,2\,1\,3\,1\,1\,1\,2\,1\,3\,0\,0\,0\,2\,1\,2\,3\,4\,3\,1\,0\,2\,0\,2\,4$\\
$9$ & $26$ & $16$ & $0$ & $31$ & $3\,3\,4\,2\,3\,2\,1\,3\,2\,0\,3\,1\,4\,2\,2\,0\,0\,3\,1\,1\,0\,1\,3\,3\,4\,4\,4\,4\,3\,3\,0\,4$\\
$10$ & $27$ & $16$ & $0$ & $31$ & $4\,3\,4\,0\,1\,2\,0\,4\,2\,1\,0\,4\,0\,0\,0\,0\,4\,0\,2\,0\,3\,4\,0\,4\,2\,0\,4\,0\,4\,4\,1\,2$\\
$11$ & $28$ & $16$ & $0$ & $30$ & $1\,4\,0\,3\,4\,2\,4\,4\,3\,1\,3\,0\,4\,4\,3\,1\,4\,2\,3\,4\,4\,4\,1\,3\,4\,4\,1\,3\,0\,0\,3$\\
$12$ & $29$ & $16$ & $0$ & $30$ & $2\,2\,2\,0\,3\,1\,1\,0\,4\,3\,2\,3\,1\,0\,4\,0\,1\,1\,1\,3\,3\,3\,0\,0\,3\,3\,1\,2\,4\,2\,2$\\
$13$ & $30$ & $16$ & $0$ & $29$ & $2\,1\,3\,3\,0\,1\,2\,2\,2\,2\,2\,4\,3\,3\,0\,4\,4\,2\,0\,3\,3\,1\,1\,4\,3\,1\,4\,0\,3\,1$\\
$14$ & $26$ & $16$ & $0$ & $24$ & $4\,0\,2\,2\,0\,4\,2\,4\,4\,4\,1\,2\,3\,4\,0\,4\,0\,0\,4\,1\,4\,2\,3\,0\,3$\\
$15$ & $27$ & $16$ & $0$ & $22$ & $4\,1\,1\,2\,1\,1\,3\,0\,3\,1\,4\,4\,1\,3\,2\,1\,2\,4\,1\,2\,4\,0\,1$\\
\bottomrule
\end{tabular}
\end{center}
\begin{center}\scriptsize
$\sigma_{4}(w)=\displaystyle\sum_{j=0}^{15}\frac{N_j(t)}{t^{a_j}(t+4)^{b_j}(t^2+4t+1)^{e_j}}\,w^j$
\\[4pt]
\begin{tabular}{@{}rrrrrl@{}}
\toprule
$j$ & $a_j$ & $b_j$ & $e_j$ & $\deg N_j$ & coefficients of $N_j$, from $t^{\deg N_j}$ down to $t^0$\\
\midrule
$0$ & $19$ & $8$ & $0$ & $25$ & $1\,1\,4\,2\,4\,4\,4\,1\,0\,2\,1\,1\,4\,0\,2\,2\,0\,4\,2\,0\,2\,0\,2\,4\,4\,1$\\
$1$ & $20$ & $9$ & $0$ & $27$ & $2\,1\,4\,2\,1\,3\,4\,4\,4\,3\,2\,4\,3\,0\,2\,3\,3\,4\,0\,4\,0\,2\,3\,0\,4\,3\,0\,3$\\
$2$ & $21$ & $9$ & $0$ & $28$ & $1\,1\,0\,2\,0\,0\,0\,4\,2\,1\,1\,3\,3\,0\,0\,0\,2\,0\,1\,3\,4\,2\,0\,4\,3\,3\,4\,3\,2$\\
$3$ & $22$ & $10$ & $0$ & $28$ & $2\,0\,1\,4\,0\,0\,2\,1\,3\,3\,1\,1\,3\,4\,4\,1\,0\,2\,4\,2\,0\,1\,0\,0\,0\,0\,3\,2\,3$\\
$4$ & $23$ & $11$ & $0$ & $30$ & $3\,3\,0\,0\,1\,2\,4\,2\,1\,2\,0\,3\,2\,4\,0\,4\,2\,3\,1\,0\,3\,2\,1\,3\,1\,2\,2\,4\,1\,3\,4$\\
$5$ & $29$ & $13$ & $0$ & $36$ & $3\,2\,0\,1\,2\,1\,0\,3\,2\,2\,3\,1\,4\,0\,2\,2\,3\,2\,1\,4\,0\,2\,4\,4\,1\,1\,3\,3\,4\,3\,3\,0\,0\,2\,4\,4\,1$\\
$6$ & $30$ & $14$ & $0$ & $38$ & $3\,4\,0\,4\,2\,0\,0\,1\,4\,0\,1\,4\,4\,2\,1\,2\,1\,4\,1\,0\,1\,2\,3\,1\,0\,1\,3\,0\,1\,3\,1\,3\,4\,2\,0\,4\,3\,0\,3$\\
$7$ & $31$ & $14$ & $0$ & $37$ & $2\,1\,2\,4\,2\,3\,0\,0\,3\,4\,0\,0\,2\,4\,3\,3\,1\,0\,3\,4\,4\,1\,0\,2\,4\,2\,3\,3\,4\,0\,3\,3\,1\,3\,2\,4\,4\,2$\\
$8$ & $32$ & $14$ & $0$ & $36$ & $2\,0\,2\,3\,3\,4\,2\,0\,4\,4\,2\,1\,1\,0\,4\,1\,2\,3\,1\,3\,4\,3\,0\,4\,3\,2\,4\,1\,2\,0\,2\,0\,0\,2\,3\,4\,1$\\
$9$ & $28$ & $14$ & $0$ & $31$ & $2\,2\,1\,2\,0\,3\,1\,3\,1\,0\,2\,4\,1\,0\,4\,1\,4\,1\,2\,3\,1\,0\,0\,4\,2\,1\,1\,3\,3\,2\,1\,1$\\
$10$ & $29$ & $14$ & $0$ & $31$ & $1\,3\,3\,1\,2\,0\,2\,4\,1\,4\,4\,3\,4\,1\,0\,4\,1\,4\,2\,2\,3\,2\,1\,0\,0\,2\,3\,2\,4\,0\,2\,3$\\
$11$ & $30$ & $14$ & $0$ & $31$ & $3\,4\,3\,1\,0\,4\,4\,4\,3\,4\,0\,3\,1\,4\,1\,0\,2\,1\,1\,0\,2\,1\,2\,3\,2\,2\,2\,1\,0\,0\,2\,2$\\
$12$ & $31$ & $14$ & $0$ & $30$ & $2\,2\,3\,0\,1\,0\,1\,1\,3\,2\,0\,1\,3\,2\,2\,1\,2\,3\,1\,1\,0\,3\,3\,4\,2\,3\,4\,1\,4\,1\,3$\\
$13$ & $32$ & $14$ & $0$ & $30$ & $4\,0\,2\,4\,2\,2\,3\,0\,3\,1\,1\,3\,2\,2\,2\,3\,1\,0\,1\,3\,3\,0\,2\,4\,3\,0\,3\,2\,2\,1\,4$\\
$14$ & $28$ & $14$ & $0$ & $24$ & $4\,1\,3\,2\,2\,1\,3\,2\,2\,2\,3\,3\,0\,1\,2\,4\,2\,0\,4\,0\,4\,3\,2\,2\,2$\\
$15$ & $29$ & $14$ & $0$ & $24$ & $3\,0\,2\,3\,1\,3\,3\,0\,1\,1\,4\,4\,4\,2\,3\,2\,2\,0\,0\,3\,2\,0\,1\,4\,4$\\
\bottomrule
\end{tabular}
\end{center}


\begin{thebibliography}{99}

\bibitem{DvH}
M.~Derickx and M.~van Hoeij,
\emph{Gonality of the modular curve $X_1(N)$},
J.~Algebra \textbf{417} (2014), 52--71.

\bibitem{DT}
M.~Derickx and K.~Terao,
\emph{Computing class groups and gonalities of algebraic curves over finite fields},
arXiv:2602.17417.

\bibitem{DerickxMSc}
M.~Derickx,
\emph{Torsion points on elliptic curves and gonalities of modular curves},
Master's thesis, Mathematisch Instituut, Universiteit Leiden, 2012,
\url{https://www.maartenderickx.nl/articles/MasterScriptie.pdf}.

\bibitem{DEvMZB}
M.~Derickx, A.~Etropolski, M.~van Hoeij, J.~S.~Morrow and D.~Zureick-Brown,
\emph{Sporadic cubic torsion},
Algebra Number Theory \textbf{15} (2021), no.~7, 1837--1864.

\bibitem{DM}
P.~Deligne and D.~Mumford,
\emph{The irreducibility of the space of curves of given genus},
Publ.\ Math.\ Inst.\ Hautes \'Etudes Sci.\ \textbf{36} (1969), 75--109.

\bibitem{DS}
M.~Derickx and A.~V.~Sutherland,
\emph{Torsion subgroups of elliptic curves over quintic and sextic number
fields}, Proc.\ Amer.\ Math.\ Soc.\ \textbf{145} (2017), no.~10, 4233--4245.

\bibitem{Conrad}
B.~Conrad,
\emph{Arithmetic moduli of generalized elliptic curves},
J.\ Inst.\ Math.\ Jussieu \textbf{6} (2007), no.~2, 209--278.

\bibitem{DR}
P.~Deligne and M.~Rapoport,
\emph{Les sch\'emas de modules de courbes elliptiques},
in: Modular Functions of One Variable~II, Lecture Notes in Math.\ \textbf{349},
Springer, Berlin, 1973, pp.~143--316.

\bibitem{Deuring}
M.~Deuring,
\emph{Reduktion algebraischer Funktionenk\"orper nach Primdivisoren des
Konstantenk\"orpers},
Math.\ Z.\ \textbf{47} (1942), 643--654.

\bibitem{Frey}
G.~Frey,
\emph{Curves with infinitely many points of fixed degree},
Israel J.\ Math.\ \textbf{85} (1994), 79--83.

\bibitem{Hess}
F.~Hess,
\emph{Computing Riemann--Roch spaces in algebraic function fields and related topics},
J.~Symbolic Comput.\ \textbf{33} (2002), 425--445.

\bibitem{Husemoller}
D.~Husem\"oller,
\emph{Elliptic Curves}, 2nd ed., Graduate Texts in Mathematics \textbf{111},
Springer, New York, 2004.

\bibitem{KM}
N.~M.~Katz and B.~Mazur,
\emph{Arithmetic Moduli of Elliptic Curves},
Annals of Mathematics Studies \textbf{108}, Princeton University Press,
Princeton, NJ, 1985.

\bibitem{Kubert}
D.~S.~Kubert,
\emph{Universal bounds on the torsion of elliptic curves},
Proc.\ London Math.\ Soc.\ (3) \textbf{33} (1976), 193--237.

\bibitem{Magma}
W.~Bosma, J.~Cannon and C.~Playoust,
\emph{The Magma algebra system. I. The user language},
J.~Symbolic Comput.\ \textbf{24} (1997), 235--265.

\bibitem{vHX142}
M.~van Hoeij,
\emph{Magma file \texttt{X1\_42.m}: lower bound for the gonality of $X_1(42)$},
\url{https://www.math.fsu.edu/~hoeij/files/X1N/LowerBoundsGonality/X1_42_57/X1_42.m}.

\bibitem{Isolated}
A.~Bourdon,
\emph{On the finiteness of isolated $j$-invariants for $X_1(N)$},
arXiv:2608.19354.

\bibitem{Najman}
F.~Najman,
\emph{Classification of torsion of elliptic curves over quintic fields},
arXiv:2609.12846.

\bibitem{NV}
F.~Najman and M.~Varivoda,
\emph{Torsion groups of elliptic curves that appear infinitely often over
septic, octic and nonic fields},
arXiv:2602.03513.

\bibitem{SuthOpt}
A.~V.~Sutherland, \emph{Optimized equations for $X_1(N)$} (joint work with
M.~van Hoeij), \url{https://math.mit.edu/~drew/X1_optcurves.html}, updated
January~10, 2014.

\bibitem{SuthMC}
A.~V.~Sutherland,
\emph{Constructing elliptic curves over finite fields with prescribed
torsion}, Math.\ Comp.\ \textbf{81} (2012), 1131--1147.

\bibitem{Sutherland}
A.~V.~Sutherland,
\emph{Defining equations for $X_1(N)$},
\url{https://math.mit.edu/~drew/X1_altcurves.html}; the model used here is
the file \url{https://math.mit.edu/~drew/FFFc42.txt}.

\end{thebibliography}
\end{document}